\documentclass{mcma_v02}
\usepackage{amssymb, url}
\usepackage{lscape}
\usepackage{color}
\usepackage{graphicx}
\usepackage{comment}
\usepackage{ascmac} 
\usepackage{mathrsfs}
\usepackage{mathptmx}      
\title{A search for extensible low-WAFOM point sets}
\headlinetitle{A search for extensible low-WAFOM point sets}

\authorone{Shin Harase}
\addressone{College of Science and Engineering, Ritsumeikan University, 1-1-1 Nojihigashi, Kusatsu, Shiga, 525-8577}
\countryone{Japan}
\emailone{harase@fc.ritsumei.ac.jp}

\researchsupported{The author was partially supported by the Grants-in-Aid for Young Scientists (B) $\#$26730015, for JSPS Fellows 24$\cdot$7985, 
for Scientific Research (B) $\#$26310211, and for Challenging Exploratory Research $\#$15K13460 from the Japan Society for the Promotion of Scientific Research. 
}

\abstract{Matsumoto, Saito, and Matoba recently proposed the Walsh figure of merit (WAFOM),
which is a computable criterion for quasi-Monte Carlo point sets using digital nets. 
Several algorithms have been proposed for finding low-WAFOM point sets.  
In the existing algorithms, the number of points is fixed in advance, 
but extensible point sets are preferred in some applications. 
In this paper, we propose a random search algorithm for extensible low-WAFOM point sets. 
For this, we introduce a method that uses lookup tables to compute WAFOM faster. 
Numerical results show that our extensible low-WAFOM point sets are comparable with Niederreiter--Xing sequences for some low-dimensional and smooth test functions.}

                                        

\firstpage{1}                                %
\pubyear{}                                   %
\volume{}                                    %
\issue{}                                     %
\doi{}                                       %
\communicated{}                              %
\received{June 20, 2016}                                  %
\revised{}                                   %

\begin{document}

\section{Introduction}\label{sec:intro}
For a Riemann integrable function $f:[0,1)^s \to \mathbb{R}$,
we consider the integral $\int_{[0,1)^s} f(\mathbf{x}) \textrm{d} \mathbf{x}$ 
and its approximation by quasi-Monte Carlo integration,  
\begin{eqnarray} \label{eqn:monte carlo}
 \int_{[0,1)^s} f(\mathbf{x}) \textrm{d} \mathbf{x} \approx \frac{1}{N}\sum_{k = 0}^{N-1} f(\mathbf{x}_k),
\end{eqnarray}
where the point set $P :=\{ \mathbf{x}_0, \ldots, \mathbf{x}_{N-1} \} \subset [0,1)^s$ is deterministically chosen 
and $N$ denotes the cardinality of $P$. 
Several criteria for $P$ have been proposed, such as the $t$-value of a digital net related to the star discrepancy \cite{MR3085113,MR918037,MR1172997,Pirsic2001827} and dyadic diaphony \cite{MR1450924}. 
Recently, Matsumoto, Saito, and Matoba \cite{MSM} proposed the {\it Walsh figure of merit} ({\it WAFOM}), 
which is a computable criterion for digital nets $P$. 
This criterion is based on Dick's error bound 
to ensure higher-order convergence than classical quasi-Monte Carlo sets 
for very smooth functions \cite{MR2346374,MR2391005,MR2683394}
and requires $O(nsN)$ arithmetic operations, where $n$ is the number of digits of precision. 
Thus, WAFOM is expected to be applied to a wide range of design and assessment for $P$. 
Generally speaking, random search is easier than a mathematical construction. 
In coding theory, low-density parity-check (LDPC) codes, which are based on random search, have been significantly successful. 
Matsumoto et al.\ \cite{MSM} proposed {\it sequential generators}, which are a special subclass of digital nets 
and only require $O(nN)$ steps for each computation of WAFOM, and obtained low-WAFOM point sets by random search in this class. 
Harase \cite{MR3314457} proposed a search algorithm for point sets whose $t$-value and WAFOM are both small 
by using random linear scrambling \cite{MR1659004} so as to be effective for a wide range of function classes. 
However, in these frameworks, the number of points $N$ has to be fixed in advance.
Extensible point sets, which are point sets with the property that the number of points may be increased 
while retaining the existing points, are preferred in some applications 
because $N$ can be increased without discarding previous function evaluations. 

The aim of this paper is to give a search algorithm for extensible low-WAFOM point sets 
based on general digital nets in the same spirit as those mentioned above. 
For this, we introduce a method that uses lookup tables to compute WAFOM values, 
which requires $O(sN)$ steps. 
Numerical experiments confirm that our extensible low-WAFOM point sets are comparable with 
(or even slightly superior to) Niederreiter--Xing sequences \cite{MR1358190} 
for some low-dimensional smooth functions.

The rest of this paper is organized as follows. In Section~\ref{sec:notation}, 
we briefly recall the notation of digital nets and WAFOM to be used in later sections. 
Section~\ref{sec:extensible low-WAFOM point sets} is devoted to our main results: 
an acceleration method that uses lookup tables to compute WAFOM
and a search algorithm for obtaining extensible low-WAFOM point sets. 
In Section~\ref{sec:Genz}, we report comparisons between our extensible low-WAFOM point sets and other quasi-Monte Carlo point sets (e.g., Sobol' \cite{MR0219238} and Niederreiter--Xing \cite{MR1358190} sequences) for the Genz test function package \cite{Genz1984, Genz1987}. 

\section{Notation} \label{sec:notation}
\subsection{Digital nets}

Let $s$ and $n$ be positive integers. Let $\mathbb{F}_2:=\{ 0,1\}$ be the two-element field, 
and let $V:= M_{s, n} (\mathbb{F}_2)$ the set of $(s \times n)$-matrices with coefficients in $\mathbb{F}_2$. 
Thus, $\mathbf{x} \in V$ denotes a matrix $\mathbf{x}:=(x_{i, j})_{1 \leq i \leq s, 1 \leq j \leq n}$ with $x_{i, j} \in \mathbb{F}_2$. 
We identify $\mathbf{x}:=(x_{i, j})_{1 \leq i \leq s, 1 \leq j \leq n}$ with an $s$-dimensional point 
$(\sum_{j=1}^{n} x_{1, j}2^{-j}, \ldots, \sum_{j=1}^{n} x_{s, j}2^{-j}) \in [0,1)^s$, which is also denoted by $(x_1, \ldots, x_s) \in [0,1)^s$.
Note that $n$ is the number of digits of precision. 

To construct $P:=\{ \mathbf{x}_0, \mathbf{x}_1,\ldots, \mathbf{x}_{2^m-1} \}$, 
we use the following construction scheme called a {\it digital net} (see \cite{MR2683394} and \cite{MR1172997} for details). 
We first select $(n \times m)$-{\it generating matrices} $C_1, \ldots, C_s \in M_{n, m}(\mathbb{F}_2)$. 
For $k = 0, 1, \ldots, 2^m-1$, let $k = \sum_{j = 0}^{m-1} k_j 2^j$ with $k_j \in \mathbb{F}_2$ 
be the expansion of $k$ in base $2$. We set $\mathbf{k} := {}^t(k_0, \ldots, k_{m-1}) \in \mathbb{F}_2^m$, where ${}^t$ represents the transpose,  
and set $\mathbf{x}_k := {}^t(C_1 \mathbf{k}, \cdots, \ldots, C_s \mathbf{k}) \in M_{s, n}(\mathbb{F}_2)$.
The point set $P:= \{ \mathbf{x}_0, \ldots, \mathbf{x}_{2^m-1} \}$ is called a {\it digital net} over $\mathbb{F}_2$. 
Throughout this paper, we assume that $P$ is a digital net, so $P \subset V$ is an $\mathbb{F}_2$-linear subspace of $V$. 
Note that the first $2^d$ points $P_d := \{ \mathbf{x}_0, \ldots, \mathbf{x}_{2^d-1} \}$ for $d \leq m$ are generated from the first $d$ columns of $C_1, \ldots, C_s$.
Note also that $P_d$ is {\it extensible}, that is, $P_d \supset P_{d-1}$. 
The quality of the point sets is determined by the generating matrices. 

\subsection{Walsh figure of merit (WAFOM)}

Matsumoto et al.\ \cite{MSM} recently proposed the Walsh figure of merit (WAFOM) 
as a computable criterion for digital nets $P$. 
Based on Dick's (Koksma--Hlawka type) inequality for integration errors \cite{MR2346374, MR2391005}, 
WAFOM has the potential to achieve higher-order convergence than $O(N^{-1})$ for function classes with very high smoothness (so-called {\it $n$-smooth functions}). 
More recently, Yoshiki \cite{Yoshiki2014} gave a more explicit error bound than Dick's by using the dyadic difference.  
Thus, throughout this paper, we adopt Yoshiki's new result and consider the same setting as in \cite{MR3314457} (with some abuse of notation). 
Following \cite{MR3314457, MatsumotoOhori2014, MSM}, we briefly recall WAFOM. 

For $\mathbf{x}=(x_{i, j})_{1 \leq i \leq s ,1 \leq j \leq n} \in M_{s, n}(\mathbb{F}_2)$, 
we define the $s$-dimensional subinterval $\mathbf{I}_{\mathbf{x}} \subset [0,1)^s$ by 
\[ \mathbf{I}_\mathbf{x} := [\sum_{j=1}^{n} x_{1, j}2^{-j}, \sum_{j=1}^{n} x_{1, j}2^{-j} + 2^{-n}) \times \cdots \times
[\sum_{j=1}^{n} x_{s, j}2^{-j}, \sum_{j=1}^{n} x_{s, j}2^{-j} + 2^{-n}).\] 
For a Riemann integrable function $f: [0,1)^s \to \mathbb{R}$, we define the $n$-digit discretization 
$f_n: V \to \mathbb{R}$ by $f_n(\mathbf{x}) := (1/{\rm Vol}(\mathbf{I}_{\mathbf{x}})) 
\int_{\mathbf{I}_{\mathbf{x}}} f(\mathbf{x}) \textrm{d} \mathbf{x}$.
This is the average value of $f$ over $\mathbf{I}_{\mathbf{x}}$. 
Under Lipschitz continuity of $f$, 
it was shown in \cite{MSM} that the discretization error between $f$ and $f_n$ on $\mathbf{I}_{\mathbf{x}}$ is 
negligible if $n$ is large enough (e.g., when $n \geq 30$). 
Thus, for $f: [0, 1)^s \to \mathbb{R}$ and large $n$, 
we may assume that $({1}/|P|) \sum_{\mathbf{x} \in P} f(\mathbf{x}) \approx (1/{|P|}) \sum_{\mathbf{x} \in P} f_n(\mathbf{x})$.

Suppose that $f$ is an $n$-smooth function (see \cite{MR2391005} and \cite[Ch.~14.6]{MR2683394} for the definition). 
Yoshiki \cite{Yoshiki2014} gave the following Koksma--Hlawka type inequality by improving Dick's inequality (\cite[Section~4.1]{MR2743889} and \cite[(3.7)]{MSM}): 
\begin{eqnarray} \label{ineq:yoshiki}
\left| \int_{[0,1)^s} f(\mathbf{x}) \textrm{d} \mathbf{x}-\frac{1}{|P|} \sum_{\mathbf{x} \in P} f_n(\mathbf{x}) \right| \leq  \sup_{0 \leq N_1, \ldots, N_s \leq n} 
|| f^{(N_1, \ldots, N_s)}||_{\infty} \cdot \mbox{WAFOM}(P),
 \end{eqnarray}
where $||f||_{\infty}$ is the infinity norm of $f$, 
and $f^{(N_1, \ldots, N_s)} := \partial^{N_1 + \cdots + N_s}f / \partial x_1^{N_1} \cdots \partial x_s^{N_s}$. 
The value $\mbox{WAFOM}(P)$ is a quality criterion for the digital net $P$ and 
is called the Walsh figure of merit (WAFOM) of $P$. 
Inequality (\ref{ineq:yoshiki}) claims that digital nets with smaller $\mbox{WAFOM}(P)$ values have better performance. 
Therefore, we want to find $P$ for which  $\mbox{WAFOM}(P)$ is small. 
This value can be obtained using the following formula:
\begin{eqnarray} \label{eqn:WAFOM}
\mbox{WAFOM}(P) = \frac{1}{|P|} \sum_{\mathbf{x} \in P} \left\{ \prod_{1 \leq i \leq s} \prod_{1 \leq j \leq n} (1+(-1)^{x_{i, j}}2^{-(j+1)}) -1 \right\}. 
\end{eqnarray} 
Thus, $\mbox{WAFOM}(P)$ is computable in $O(nsN)$ arithmetic operations, where $\mathbf{x} := (x_{i, j}) \in P$ and $N = |P|$. 
We refer the reader to \cite{MSM} and \cite[Section~2]{MatsumotoOhori2014} for details. 
In the next section, we search for extensible low-WAFOM point sets by means of random search.

\section{A search for extensible low-WAFOM point sets} \label{sec:extensible low-WAFOM point sets}
\subsection{Acceleration using lookup tables} \label{subsec:table}

To search for a low-WAFOM point set $P$, 
it is crucial to compute (\ref{eqn:WAFOM}) as fast as possible. 
Matsumoto et al.~\cite{MSM} restricted the search-space to 
a special subclass of digital nets, called {\it sequential generators}, and reduced the complexity to 
$O(nN)$ arithmetic operations for computing $\mbox{WAFOM}(P)$ for each $P$.  
However, for general digital nets, we need $O(nsN)$ steps for the na\"ive method in (\ref{eqn:WAFOM}).
In this subsection, we propose the use of lookup tables as another direction of acceleration.  

For simplicity, we assume 
without loss of generality that 
$q$ is a positive integer that divides $n$, e.g., $n= 30$ and $q = 3$. 
We consider the $i$-th row vector $\mathbf{x}^{(i)} := (x_{i, 1}, \ldots, x_{i, n}) \in \mathbb{F}_2^n$ of $\mathbf{x} \in P$.   
We split $\mathbf{x}^{(i)}$ into $q$ equal segments, and set $l := n/q$.  
Thus, $\mathbf{x}^{(i)}$ is decomposed into $(\mathbf{d}_{1}^{(i)}, \ldots, \mathbf{d}_{q}^{(i)})$, 
where $\mathbf{d}_{c}^{(i)} := (x_{i, (c-1)l+1}, \ldots, x_{i, cl})$ for each $c = 1, \ldots, q$, that is,
\begin{eqnarray*}
\mathbf{x}^{(i)} = \underset{\mathbf{d}_{1}^{(i)}}{\underbrace{x_{i, 1}, \ldots, w_{i, l}}}, \underset{\mathbf{d}_{2}^{(i)}}{\underbrace{x_{i, l+1}, \ldots, x_{i, 2l}}}, \ldots, \underset{\mathbf{d}_{q}^{(i)}}{\underbrace{x_{i, n-l+1}, \ldots, x_{i, n}}}.
\end{eqnarray*}
We identify $\mathbf{d}_{c}^{(i)}$ with the $l$-bit integer $\sum_{j=1}^l x_{i, {(c-1)l + j}} 2^{l-j}$, 
so that $\mathbf{d}_{c}^{(i)}$ is viewed as an integer from $0, 1, \ldots, 2^{l-1}$.

For each $c = 1, \ldots, q$,  
we construct the following lookup tables in advance: 
\begin{eqnarray} \label{eqn:table}
{\tt table}_{c}[\mathbf{e}] \gets \prod_{1 \leq j \leq l} (1+(-1)^{e_{j}}2^{-((c-1)l + j+1)}) \quad (\mathbf{e} = 0, \ldots, 2^l-1), 
\end{eqnarray}
where each $e_j \in \mathbb{F}_2$ is given by the binary expansion $\mathbf{e} = \sum_{j = 1}^{l} e_j2^{l-j}$.  
The products $\prod_{1 \leq j \leq n} (1+(-1)^{x_{i, j}}2^{-(j+1)})$ in (\ref{eqn:WAFOM}) then reduce to 
$\prod_{1 \leq c \leq q} {\tt table}_{c}[\mathbf{d}_{c}^{(i)}]$.  
Thus, we obtain the following proposition.
\begin{proposition} \label{prop:table}
If we use the above lookup tables, then (\ref{eqn:WAFOM}) is computable in $O(sN)$ arithmetic operations, 
where the hidden constant depends on $q$. 
In other words, our method decreases the number of multiplications by a factor of $1/l$.
\end{proposition}

For increasing speed,
we recommend selecting a small $q$ such that 
the corresponding tables made by (\ref{eqn:table}) are included in cache memory. 
(We also recommend using $q$ 1-dimensional arrays rather than a 2-dimensional array as lookup tables.) 
To show the effectiveness of our approach, we conduct experiments comparing the na\"ive method and our lookup-table method. 
We set $(n, q) = (30, 3)$ and generate the first $2^{25}$ points $P$ of the Niederreiter--Xing sequences 
implemented by Pirsic \cite{MR1958872} with Gray code order \cite{Bratley:1992:ITL:146382.146385}. 
We measure the CPU time for computing WAFOM values for $s = 4, 6, 8, \ldots, 16$.
The experiments are conducted on a 64-bit Intel Core i7-3770 3.90 GHz CPU.
Our codes are implemented in C and compiled by the GCC compiler with the -O3 optimization flag on a Linux operating system. 
As shown in the second and third rows of Table~\ref{table:timing}, in each case, our method runs over $30$ times faster than the na\"ive method. 
By using the lookup tables, we also have the advantages that the arithmetic operations and the conditional branches in (\ref{eqn:WAFOM}) are avoided.

Ohori and Yoshiki \cite{OY} recently gave a fast and simple method for computing a good approximation of WAFOM that was originally proposed in \cite{MSM}. 
In summary, ${\rm WAFOM}(P)$ is well-approximated by the QMC-error of the function $f(\mathbf{x}) = \exp(-2\sum_{i = 1}^{s} x_i)$ 
(see Remark~2 in \cite{MatsumotoOhori2014}). 
We implement their method on the above platform, 
and record the CPU times in the last row of Table~\ref{table:timing}.
Our lookup table method is still faster than the approximation method when $s$ is small. 
This is possibly because the calculation of $\exp(x)$ is much slower than multiplication on some platforms. 
We note that these timings strongly depend on the CPU  (see, for example, Chapter~3.7 in \cite{Patterson:2013:COD:2568134}). 
For randomized digital nets, no corresponding approximation method is known. See Remark~\ref{remark:RQMC} for details.

\begin{table}[h]
\begin{center}
\caption{CPU time (sec) taken to compute WAFOM values for $2^{25}$-point Niederreiter--Xing sequences ($n = 30$).}
\begin{tabular}{c||c|c|c|c|c|c|c} \hline
$s$ & $4$ & $6$ & $8$ & $10$ & $12$ & $14$ & $16$\\ \hline 
{Na\"ive} & $11.773$ & $19.307$ & $26.651$ & $33.039$ & $39.527$ & $46.119$ & $52.639$ \\ \hline 
{Table ($q = 3$)} & $0.367$ & $0.505$ & $0.654$ & $0.795$ & $0.972$ & $1.112$ & $1.249$ \\ \hline
{Approximation} & $1.111$ & $1.197$ & $1.282$ & $1.384$ & $1.479$ & $1.565$ & $1.676$ \\ \hline
\end{tabular}
\label{table:timing}
\end{center}
\end{table}

\subsection{A heuristic algorithm for searching for extensible low-WAFOM point sets} \label{subsec:searching}

To obtain an extensible point set $P_d := \{ \mathbf{x}_0, \ldots, \mathbf{x}_{2^d-1} \} \supset P_{d-1}$, 
we consider a search algorithm for determining $(n \times m)$-generating matrices $C_1, \ldots, C_s \in M_{n, m}(\mathbb{F}_2)$ 
for which $\mbox{WAFOM}(P_d)$ is small
for each $d = 1, \ldots, m$. Let us decompose $C_i$ into column vectors so that $C_i = ({\overline{\mathbf{c}}}_1^{(i)}, \ldots, {\overline{\mathbf{c}}}_m^{(i)})$. 
Our approach is as follows. 
\begin{enumerate}
\item[(1)] We generate the first column vectors ${\overline{\mathbf{c}}}_1^{(1)}, \ldots, {\overline{\mathbf{c}}}_{1}^{(s)}$ 
at random $M$ times, 
and select the first column vectors with the smallest $\mbox{WAFOM}(P_1)$ values.
\item[(2)] For $d = 2, \ldots, m$, if the first $(d-1)$ column vectors
${\overline{\mathbf{c}}}_1^{(i)}, \ldots, {\overline{\mathbf{c}}}_{d-1}^{(i)} \in \mathbb{F}_2^n$ ($i = 1, \ldots, s$) have been determined, 
then we generate the $d$-th column vectors ${\overline{\mathbf{c}}}_d^{(1)}, \ldots, {\overline{\mathbf{c}}}_d^{(s)}$ at random $M$ times, 
and select the $d$-th column vectors with the smallest $\mbox{WAFOM}(P_d)$ values. 
\end{enumerate}
To ensure the uniformity of each 1-dimensional projection, 
we generate ${\overline{\mathbf{c}}}_d^{(1)}, \ldots, {\overline{\mathbf{c}}}_d^{(s)}$ so that 
the upper $(d \times d)$-submatrix of 
$({\overline{\mathbf{c}}}_1^{(i)}, \ldots, {\overline{\mathbf{c}}}_{d}^{(i)})$ is regular for each $i = 1, \ldots, s$. 
This property ensures that each 1-dimensional projection is optimal form the view point of the $t$-value of the $(t, m, s)$-net  
(that is, the $t$-value is $0$). See \cite {MR2683394,MR1172997} for details. 

\begin{figure} 
\centering
  \includegraphics[width=11cm]{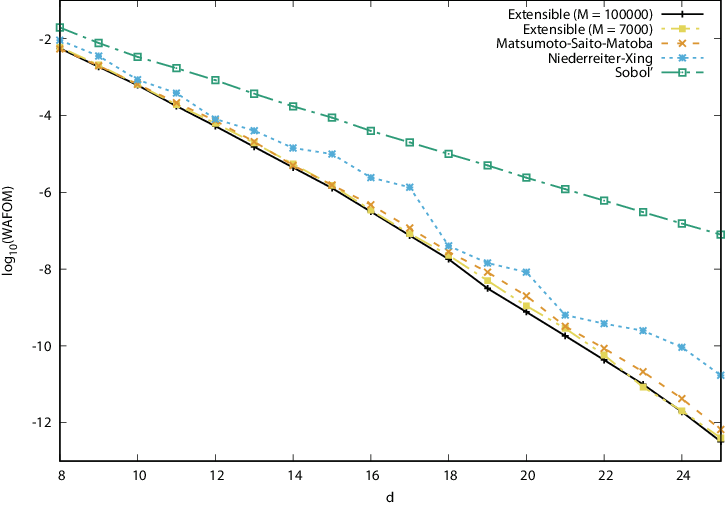}
  \caption{WAFOM (on a $\log_{10}$ scale) for our extensible point sets (with $M = 7000$ and $M = 100000$) and point sets from existing methods.}
\label{fig:WAFOM}
\end{figure}

We compare our methods with the existing methods including the sequential generators from \cite{MSM}. 
From now on, we set $n = 32$. We search for sequential generators in the same way as Matsumoto et al.\ \cite{MSM}. 
For this, we conduct $7000$ random searches for $(n, s) = (32, 5)$. 
Conversely, for our methods, we set $(n, m, s) = (32, 25, 5)$ and search for column vectors using $M = 7000$ trials and $M = 100000$ trials. 
We also compute WAFOM values for the Sobol' \cite{MR2429482} and Niederreiter--Xing \cite{MR1958872} sequences. 
Figure~\ref{fig:WAFOM} shows the WAFOM values for these point sets for $d = 8, \ldots, 25$. 
Note that the sequential generators do not possess extensibility.
For $M = 7000$ trials, 
our extensible low-WAFOM point sets are slightly better than the point sets constructed from sequential generators for large $d$. 
This is possibly because Matsumoto--Saito--Matoba sequential generators have a somewhat restricted parameter space.

\begin{remark} \label{remark:RQMC}
In terms of the Walsh coefficients of $f$, 
Goda, Ohori, Suzuki, and Yoshiki \cite{GOSY} proposed a quality criterion $\mathscr{W}(P)$ for randomized quasi-Monte Carlo integration 
using digitally shifted digital nets. This criterion is called the {\it Walsh figure of merit for root mean square error}, 
which satisfies a Koksma--Hlawka type inequality on the root mean square error. The computable formula for $\mathscr{W}(P)$ is given by 
\[ \mathscr{W}(P) := \sqrt{ \frac{1}{|P|} \sum_{\mathbf{x} \in P} \left\{ \prod_{1 \leq i \leq s} \prod_{1 \leq j \leq n} (1+(-1)^{x_{i, j}}2^{-2(j+1)}) -1 \right\}}. \]
Our search algorithm is also applicable here by replacing $2^{-(j+1)}$ in (\ref{eqn:WAFOM}) with $2^{-2(j+1)}$. 
At present, an efficient approximation of $\mathscr{W}(P)$, such as that mentioned in 
Section~\ref{subsec:table}, is not known. 
Thus, our lookup table method seems to be a good option in the case of randomization by digital shifts as well. 
\end{remark}

\begin{remark} \label{remark:suzuki}
Suzuki \cite{MR3183338} gave an explicit mathematical construction of low-WAFOM point sets 
using Niederreiter-Xing sequences and Dick's interlacing construction \cite{MR2346374, MR2391005} 
for a fixed size $N = 2^m$. 
A mathematical construction of extensible low-WAFOM point sets is an open problem.
\end{remark}

\section{Numerical experiments} \label{sec:Genz}

We evaluate the following four methods: 
\begin{enumerate}
\item[(a)] Extensible low-WAFOM point sets for $M=100000$ (using the procedure in the previous section);
\item[(b)] Matsumoto--Saito--Matoba sequential generators \cite{MSM} (using the procedure in the previous section);
\item[(c)] Niederreiter--Xing sequences \cite{MR1358190} implemented by Pirsic \cite{MR1958872};
\item[(d)] Sobol' sequences with good two-dimensional projections \cite{MR2429482}.
\end{enumerate}
For point sets (a)--(d), 
the WAFOM values are plotted in Figure~\ref{fig:WAFOM} in the previous section. 
We use the following six different types of test function from the Genz package \cite{Genz1984, Genz1987} defined over $[0,1)^s$:
\begin{eqnarray*}
\begin{array}{ll}
\mbox{Oscillatory:} & f_1(\mathbf{x})  =  \cos (2 \pi u_1 + \sum_{i = 1}^{s} a_i x_i),\\
\mbox{Product Peak:} & f_2(\mathbf{x})  =  \prod_{i = 1}^{s} [1/{(a_i^{-2} + (x_i - u_i)^2)]},\\
\mbox{Corner Peak:} & f_3(\mathbf{x})  =  (1 + \sum_{i = 1}^{s} a_i x_i)^{-(s+1)},\\
\mbox{Gaussian:} & f_4(\mathbf{x})  =  \exp (-\sum_{i=1}^s a_i^2 (x_i - u_i)^2),\\
\mbox{Continuous:} & f_5(\mathbf{x})  =  \exp (-\sum_{i=1}^s a_i |x_i - u_i|),\\
\mbox{Discontinuous:} & f_6(\mathbf{x})  =  \left\{
\begin{array}{ll}
0, & \mbox{if $x_1 > u_1$ or $x_2 > u_2$},\\
\exp (\sum_{i = 1}^s a_i x_i), & \mbox{otherwise.}
\end{array}
\right.
\end{array}
\end{eqnarray*}
We can obtain the exact value $I(f_j) := \int_{[0,1)^s} f_j \textrm{d} \mathbf{x}$ analytically, 
so such families have been used as test functions \cite{MR1417864,MR1958872, mSLO94a} 
and analyzed from a theoretical point of view \cite{MR1963917}. 
The different test integrands are obtained by changing $\mathbf{a} = (a_1, \ldots, a_s)$ 
and $\mathbf{u} = (u_1, \ldots, u_s)$. 
The parameter $\mathbf{a}$ affects the degree of difficulty, and 
the parameter $\mathbf{u}$ acts as a shift parameter. 
We generate $\mathbf{a}$ and $\mathbf{u}$ as uniform random vectors in $[0, 1]^s$, 
and $\mathbf{a}$ is renormalized to satisfy
\begin{eqnarray} \label{eqn:Genz condition}
\sum_{i = 1}^s a_i = h_j,
\end{eqnarray}
where $h_j$ depends on the family $f_j$. 
This condition determines the difficulty of integration, 
as the difficulty of computing the integral $I(f_j)$ increases when $h_j$ increases. 
In this way, we obtain a function $f_j$ for each test family in dimension $s$ and make 20 quantitative examples of $f_j$ by changing $\mathbf{a}$ and $\mathbf{u}$. 
For any sample size $|P|=2^d$ and any family $f_j$, 
we compute the median of the relative errors (on a $\log_{10}$ scale)
\[ \log_{10} \frac{|I(f_j) - I_N(f_j)|}{|I(f_j)|}, \]
where $N := |P|$, $I_N(f_j) := ({1}/{|P|}) \sum_{\mathbf{x} \in P} f_j(\mathbf{x})$, and $I(f_j) := \int_{[0,1)^s} f_j \textrm{d} \mathbf{x}$.
As suggested in Remark~2.2 in \cite{MSM}, we translate $P$ by adding $(2^{-n-1}, \ldots, 2^{-n-1})$ to obtain better performance for numerical integration.

\begin{figure}
\includegraphics[width=8cm]{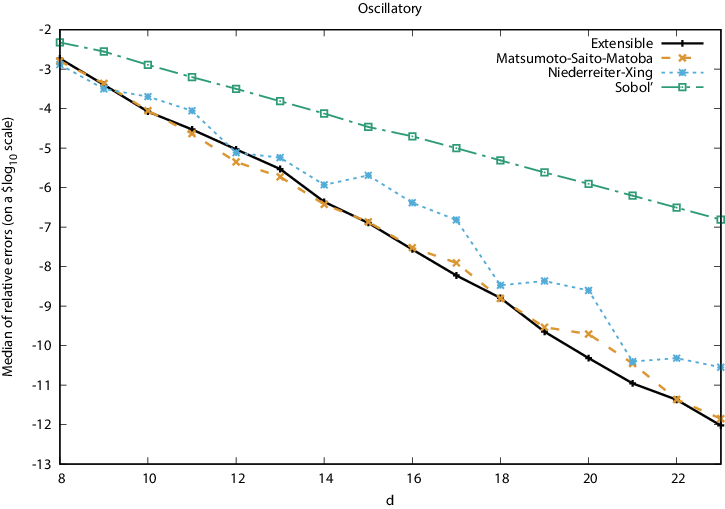}
\includegraphics[width=8cm]{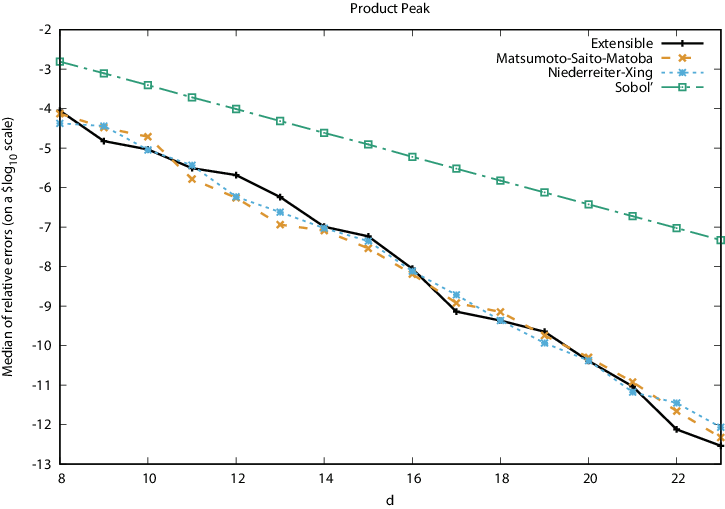}
\includegraphics[width=8cm]{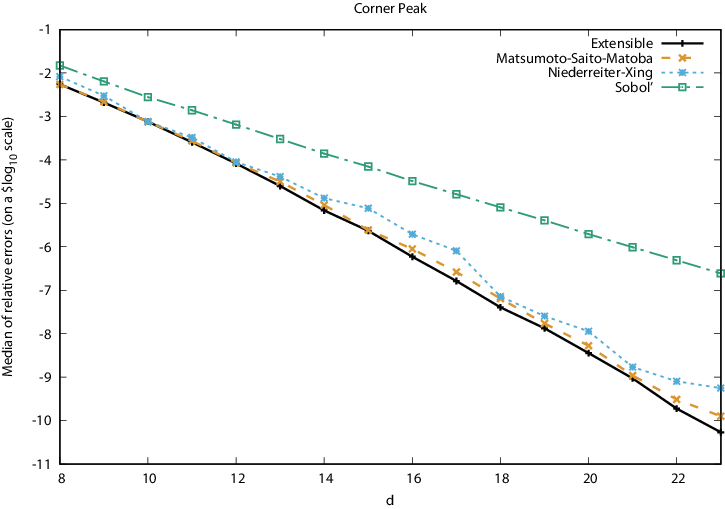}
\includegraphics[width=8cm]{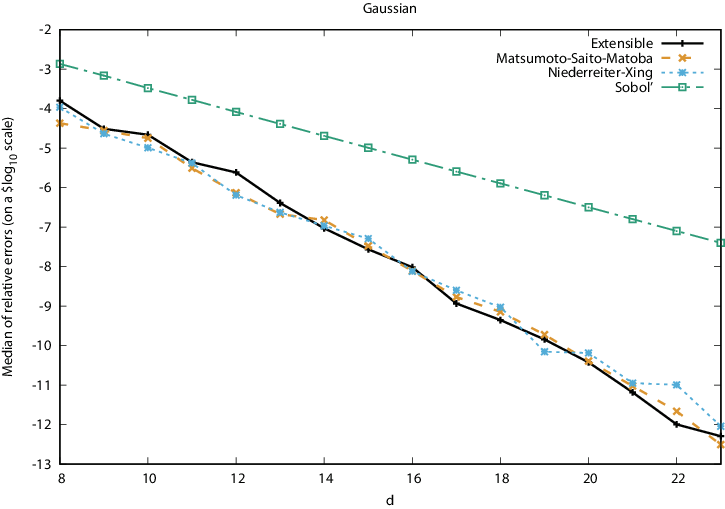}
\includegraphics[width=8cm]{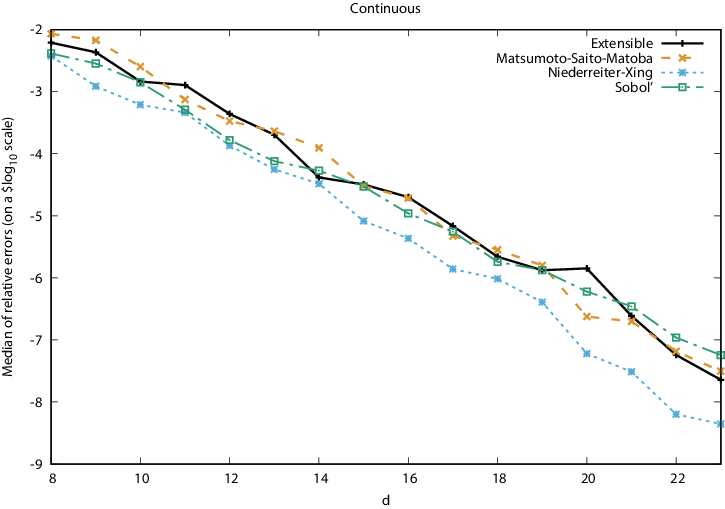}
\includegraphics[width=8cm]{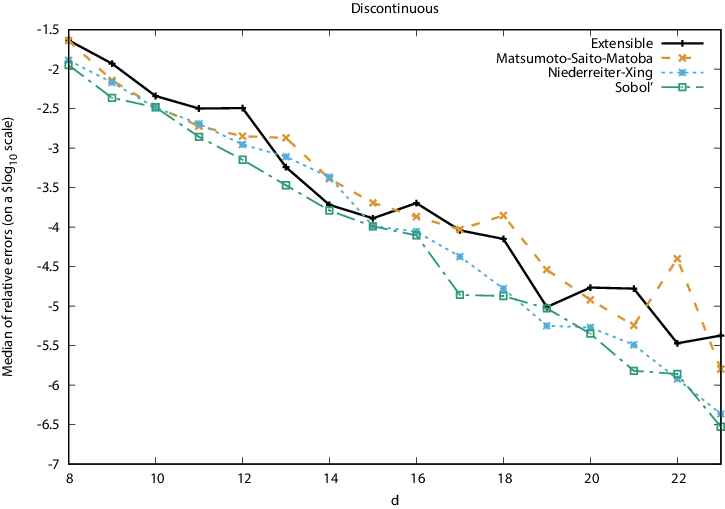}
  \caption{
The medians of the relative errors (on a $\log_{10}$ scale) for 5-dimensional Genz test functions;
the horizontal axis shows $d$, and the vertical axis the medians of the relative errors (on a $\log_{10}$ scale).}
  \label{fig:genz5}
\end{figure}

We select the parameters $s = 5$ and $(h_1, h_2, h_3, h_4, h_5, h_6) = (4.5, 3.625, 0.925, 3.515, 10.2, 2.15)$, 
which are the same settings as in \cite{MR3314457}. 
(These values are half of the values used in the case of $10$-dimensional test functions in \cite[p.~284]{MR1768951}.) 
Figure~\ref{fig:genz5} gives a summary of the median of the relative errors for $d = 8, \ldots, 23$. 
Point sets (a) and (b) are superior to the Sobol' and Niederreiter--Xing sequences 
for $f_1$ and $f_3$, and are comparable with the Niederreiter--Xing sequence for $f_2$ and $f_4$. 
However, for $f_5$ (continuous but not differentiable functions), 
the Niederreiter--Xing sequence is superior to the WAFOM-based methods.
For $f_6$, the Sobol' sequence seems to be best. 
Recall that the $t$-value of the $(t, m, s)$-net \cite{MR1172997}, which is a well-established figure of merit for a digital net.  
The Niederreiter--Xing and Sobol' sequences are optimized in terms of the $t$-values.  
For non-smooth functions, these point sets seem to be more effective than our simple low-WAFOM point sets, 
which do not consider the $t$-values (except for 1-dimensional projections). See, for example, Table~1 in \cite{MR3314457}. 
This tendency coincides with the experimental results in \cite{MatsumotoOhori2014}. 
Harase \cite{MR3314457} gave a search algorithm for quasi-Monte Carlo point sets with both small WAFOM values and $t$-values 
using random linear scrambling \cite{MR1659004} to 
improve the rates of convergence for smooth functions while being robust for non-smooth functions. 
However, his point sets do not have the property of extensibility. 

When the dimension $s$ is high, 
our extensible low-WAFOM points are inferior to the Niederreter--Xing sequences with the exceptions of $f_1$ and $f_3$. 
Although we conducted experiments with $s = 10$, 
it seems to be difficult to obtain good point sets by simple random search for high dimensions. 
(Theoretically, the smallest value of $\log({\rm WAFOM}(P_d))$ is $O(-d^2/s)$ for $P_d$ with $|P_d| = 2^d$, 
so the rate of convergence tends to be worse when $s$ is large, see \cite{MR3145585,MR3429168,MR3296075}.)
A breakthrough in high-dimensional integration is to take into account importance of variates, known as {\it weights}. 
Suzuki \cite{Suzuki2015} recently studied an infinitely differentiable function space with certain weights for which WAFOM (with weights) works well, 
and showed a sufficient condition for very fast convergence (called {\it accelerating convergence with strong tractability}) 
when the weights decay sufficiently quickly. 
Research on an efficient approximation of WAFOM with weights (called {\it WAFOM with derivation sensitivity parameter}) is also reported \cite{Ohori2015,OY}. 
We refer the reader to a recent survey for details \cite{MatsumotoOhori2014}.



\section{Conclusions}
We have proposed a random search algorithm for extensible low-WAFOM point sets in terms of digital nets. 
The key to our algorithm is determining the columns of generating matrices, inductively. 
We also introduce a lookup-table method to compute WAFOM faster. 
The point sets obtained have almost the same level of accuracy as 
Matsumoto--Saito--Matoba non-extensible point sets for numerical integration, 
and are comparable with (or even slightly superior to) Niederreiter--Xing sequences for some low-dimensional functions with high smoothness. 

\section*{Acknowledgments}
In an earlier version of this paper, Mr.\ Ryuichi Ohori checked the code in C written by the author, 
in particular for the correctness of the WAFOM values in Section~\ref{sec:extensible low-WAFOM point sets}, 
and gave the author invaluable comments to improve the presentation. The author thanks him for his help.

\bibliography{wafombib-2}

\end{document}